\documentclass[12pt]{amsart}
\theoremstyle{remark}
\newtheorem*{remark}{Remark}
\theoremstyle{fact}
\newtheorem*{fact}{Fact}
\theoremstyle{plain}

\newtheorem*{theorem}{Theorem}

\theoremstyle{definition}
\newtheorem*{definition}{Definition}
\newtheorem*{Acknowledgement}{Acknowledgement}
\theoremstyle{acknowledgements}

\newcommand{\V}{\Vert}

\def\cc#1{\hbox to14pt{\hfil${\textstyle{#1}}$\hfil}}
\begin{document}
\textwidth=30cc
\baselineskip=16pt
\title{A Conditional Quasi-greedy Basis of $l_1$}
\author{S. J. Dilworth and David\ Mitra}
\address{Department of Mathematics, University of South Carolina,
         Columbia, SC 29208, U.S.A.}
\email{dilworth@math.sc.edu}\email{mitra@math.sc.edu}
\date{August 13, 2000}
\begin{abstract}
      We show that the Lindenstrauss basic sequence in $l_1$ may be used to construct a conditional
      quasi-greedy basis of $l_1$, thus answering a question of Wojtaszczyk. We further show
      that the sequence of coefficient functionals for this basis is not quasi-greedy.
\end{abstract}
\maketitle
\section{Introduction}

In what follows $\{e_i\}_{i=1}^\infty$ denotes the standard unit vector basis of $l_1$.
In~\cite{W}, the following concept was studied:

\begin{definition} Let $\frak X$ be a Banach space with dual $\frak X^*$ and let
$(x_i, x^*_i)_{i\in F}$ be a fundamental biorthogonal
system in $\frak X\times \frak X^*$ with $\inf\limits_{i\in F} \V x_i\V>0$ and
$\sup\limits_{i\in F} \V x_i^*\V<\infty$. For $m\in \mathbb{N}$, define the
operator $\mathcal{G}_m$ by
$$
  \mathcal{G}_m(x) = \sum_{i\in A} x_i^*(x) x_i,
$$
where $A\subset F$ is a set of cardinality $m$ such that $|x_i^*(x)|\ge|x_k^*(x)|$
whenever $i\in A$ and $k\notin A$ (note that the set $A$
depends on $x$ and may not be unique; nevertheless, $\mathcal{G}_m$ is well-defined).
Then
$(x_i, x_i^*)_{i\in F}$ is a {\it quasi-greedy system}
provided that the operators $\mathcal{G}_m$  satisfy:
$$
  \lim_m \mathcal{G}_m(x) =x,\quad{\rm for\ each\ }x\in\frak X.
$$
Equivalently, by \cite[Theorem~1]{W}, $(x_i, x^*_i)_{i\in F}$ is a quasi-greedy
system if there exists a constant $C$ so that for every $x\in\frak X$ and for every
$m\in\Bbb N$, we have
$$
  \V\mathcal{G}_m(x)\V\le C\V x\V.
$$
If $\{x_i\}_{i=1}^\infty$ is, in addition, basic, it will be said to be
a {\it quasi-greedy basis\/.} A sequence $\{z_i\}_{i=1}^\infty$
is {\it unconditional for constant coefficients\/} if there exist positive
numbers $c$ and $C$ such that
$$
  c\biggl\V \sum\limits_{i=1}^m z_i\biggr\V\le \biggl\V\sum\limits_{i=1}^m
            \epsilon_i z_i\biggr\V\le C\biggl\V\sum\limits_{i=1}^m z_i\biggr\V
$$
for any positive integer $m$ and for any sequence of signs $\{\epsilon_i\}_{i=1}^m$.
\end{definition}

\begin{remark}
By \cite[Proposition~2]{W} a quasi-greedy basis is necessarily
unconditional for constant coefficients.
\end{remark}

In \cite{W}, Wojtaszczyk shows that, for $1<p<\infty$, the space $l_p$ possesses a
conditional quasi-greedy basis. Since, by the Remark, a quasi-greedy basis of $l_1$
is quite close to being equivalent to $\{e_i\}_{i=1}^\infty$,
and since  $\{e_i\}_{i=1}^\infty$ is the unique normalized unconditional
basis of $l_1$ up to equivalence \cite{LP}, it is of interest to find a
conditional quasi-greedy basis of $l_1$. 
In particular, the existence of such a basis would show that
$l_1$ does not have a unique (up to equivalence)
normalized basis that is unconditional for constant coefficients. In this note, we show that indeed such a basis
exists.

 Our example is derived from the basic sequence constructed by
Lindenstrauss in \cite{L}. This is a monotone, conditional basic sequence in $l_1$ whose
closed linear span is a ${\mathcal L}_1$-space possessing no unconditional basis.
We denote this sequence by $\{x_i\}_{i=1}^\infty$ and its associated sequence of coefficient
functionals by $\{x_i^*\}_{i=1}^\infty$. They are defined as follows:
For $i\in\mathbb N$,
$$
  x_i=e_i-\textstyle{1\over2}(e_{2i+1}+e_{2i+2}).
$$
Considering $x_i^*$ as an element of $l_\infty/[x_i]^\perp$,
we write $x_i^*=y_i^*+[x_i]^\perp$ where the $y_i^*\in l_\infty$ are as defined by
Holub and Retherford (see \cite{HR}): For each $i\in\mathbb N$,
let $\alpha_i$ be the finite sequence of positive integers defined by the
following conditions:
\begin{itemize}
      \item[1)] $\alpha_i(1)=i$.
      \item[2)] $\alpha_i(j)=\alpha_i(j-1) -\bigl( [\alpha_i(j-1)/2]+1\bigr)$
                for admissible $j<i$ (that is, such that $\alpha_i(j)>0$),
                where $[k]$ denotes the greatest integer less than or equal to $k$.
\end{itemize}
Then $y_i^*$ is defined by
$$
  y_i^*=\sum_{j=1}^{|\alpha_i|}{\bigl(\textstyle{1\over2}\bigr)}^{j-1} e_{\alpha_i(j)}
$$
Thus, for example,
\begin{align*}
y_1^*&=(\cc{1},\cc{0},\ldots \,)   \\
y_2^*&=(\cc{0},\cc{1},\cc{0},\ldots \,)\\
y_3^*&=(\cc{1\over2},\cc{0},\cc{1},
                    \cc{0},\ldots \,)   \\
y_4^*&=(\cc{1\over2},\cc{0},\cc{0},\cc{1},
                          \cc{0},\ldots \,)   \\
y_5^*&=(\cc{0},\cc{1\over2},\cc{0},\cc{0},
                      \cc{1},\cc{0},\ldots \,)   \\
y_6^*&=(\cc{0},\cc{1\over2},\cc{0},\cc{0},\cc{0},
                         \cc{1},\cc{0},\ldots \,)   \\
y_7^*&=(\cc{1\over4},\cc{0},\cc{1\over2},\cc{0},
                \cc{0},\cc{0},\cc{1},\cc{0},\ldots \,)   \\
y_8^*&=(\cc{1\over4},\cc{0},\cc{1\over2},
        \cc{0},\cc{0},\cc{0},\cc{0},\cc{1},\cc{0},\cc{0},\ldots \,).   \\
\end{align*}

The properties of $\{x_i\}_{i=1}^\infty$ which we shall use in the sequel
are summarized in the following fact (see \cite{L}, \cite{LP}, and, also, \cite{S}).
\begin{fact}
The sequence $\{x_i\}_{i=1}^\infty$ satisfies:
\begin{itemize}
      \item[1)] $\{x_i\}_{i=1}^\infty$  is a monotone basic sequence {\rm(\cite[p.\ 455]{S})}.
      \item[2)] $[x_i]$ has no unconditional basis {\rm(\cite[p.\ 455]{S})}.
      \item[3)] For $n\in\mathbb N$, there exists an isomorphism
                 $T_n$ from $[x_i:1\le i\le n]$ onto $l_1^n$
                satisfying $\V T_n\V\V T^{-1}_n\V\le 2$ {\rm(\cite[Ex. 8.1]{LP})}.
\end{itemize}
\end{fact}

\section{Results}
We now construct the basis heralded in the Introduction. To do this, it suffices by \cite[Proposition~3]{W} to
construct such a basis in a space isomorphic to
$l_1$. Towards this end, define, for each $n\in\Bbb N$,
$$
  F_n=[x_i :1\le i\le n].
$$
Let
$$
  \frak X=\Bigl(\sum_{i=1}^\infty\oplus F_i\Bigr)_1.
$$
We claim that the natural basis $\{\tilde x_i\}_{i=1}^\infty$ of $\frak X$ obtained from
the $x_i$'s is the desired sequence. Indeed, it follows
from the Fact that $\{\tilde x_i\}_{i=1}^\infty$ is a monotone, conditional
basis of $\frak X$, and that $\frak X$ is isomorphic to~$l_1$.
Moreover, assuming for the moment that $\{x_i\}_{i=1}^\infty$ is quasi-greedy,
for $m\in\Bbb N$ and $y=\sum y_i\in\frak X$, we have
$$
  \V{\mathcal G_m}(y)\V\le\sup_{k\in K_m}\sum_{i=1}^\infty\V{\mathcal G}_{k(i)}(y_i)\V
                  \le C\sum_{i=1}^\infty\V y_i\V=C\V y\V,
$$
where $C$ is as in the Definition,
$K_m=\{\,k:\Bbb N\rightarrow\Bbb N\cup\{0\}:\sum\limits_{i=1}^\infty k(i)=m\,\}$, and
$\mathcal{G}_0(x)=0$ for each $x\in\frak X$ (note that the operators $\mathcal G_j$ appearing
in the above inequality are defined with respect to two different sequences).
Thus, the fact that $\{ \tilde x_i\}_{i=1}^\infty$  is quasi-greedy will be
a consequence of
the following theorem.

\begin{theorem}
The sequence $\{ x_i\}_{i=1}^\infty$ satisfies
\begin{equation}\label{eq: quasigreedy}
    3  \biggl\V\sum_{i\in S_1\cup S_2}\alpha_i x_i \biggr\V\ge  \biggl\V\sum_{i\in S_1} \alpha_i x_i\biggr\V,
\end{equation}
whenever $S_1$ and $S_2$ are disjoint finite subsets of $\Bbb N$ with
\begin{equation}\label{eq: qg}
      \min_{i\in S_1}|\alpha_i|\ge\max_{i\in S_2} |\alpha_i|.
\end{equation}
\end{theorem}

\begin{proof}
For $A\subseteq\mathbb N$ and $x\in l_1$, we denote by $P_Ax$ the vector in $l_1$
whose $j^{\rm th}$-coordinate is $x(j)$ if $j\in A$ and zero
otherwise.
Let $S_1$, $S_2$, and $\{\alpha_i\}_{i\in S_1\cup S_2}$ be as above. Set
$$
  x=\sum_{i\in S_1} \alpha_i x_i\quad{\rm and}\quad y=\sum_{i\in S_2}\alpha_i x_i,
$$
and define the sets:
\begin{align*}
      A_0&=\Bigl\{\,j\in\Bbb N:{\textstyle\sum\limits_{i\in S_1}}x_i(j)=1\,\Bigr\},      \\
      B_0&=\Bigl\{\,j\in\Bbb N:{\textstyle\sum\limits_{i\in S_1}} x_i(j)=-1/2\,\Bigr\},   \\
    \intertext{and}
      C_0&=\Bigl\{\,j\in\Bbb N:{\textstyle\sum\limits_{i\in S_1}} x_i(j)=1/2\,\Bigr\}.   \\
\end{align*}

First, we concentrate our attention on the set $B_0$. We define the sets:
\begin{align*}
      W_1&=\{\,i\in S_2:x_i(j)=1{\rm\ for\ some\ }j\in B_0\,\},                           \\
      A_1&=\{\,j\in A_0:x_i(j)=-1/2{\rm\ for\ some\ }i\in W_1\,\},                          \\
    \intertext{and}
      B_1&=\{\,j\notin A_0:x_i(j)=-1/2{\rm\ for\ some\ }i\in W_1\,\}.                      \\
\end{align*}
Finally, set
$$
  y_0=x{\rm\ and\ }  y_1=\sum_{i\in W_1} \alpha_i x_i.
$$
Note that $A_1, B_1$, $B_0$, and $C_0$ are mutually disjoint.
We also have from the triangle inequality that
\begin{equation}\label{eq: oone}
  \V P_{B_0} y_0 \V \le  \V P_{B_0} (y_0+y_1 ) \V + \V P_{B_0} y_1  \V.
\end{equation}
But, since $\V P_{B_0} y_1 \V = \V P_{A_1} y_1 \V + \V P_{B_1} y_1 \V$,
we obtain from \eqref{eq: oone} that
$$
  \V P_{B_0}( y_0+y_1) \V+\V P_{B_1}y_1\V\ge\V P_{B_0} y_0\V-\V P_{A_1}y_1\V. \leqno{(*)}
$$

Concentrating our attention now on the set $B_1$, we let
\begin{align*}
       W_2&=\{\,i\in S_2:x_i(j)=1{\rm\ for\ some\ }j\in B_1\,\},          \\
       A_2&=\{\,j\in A_0:x_i(j)=-1/2{\rm\ for\ some\ }i\in W_2\,\},         \\
\intertext{and}
       B_2&=\{\,j\notin A_0:x_i(j)=-1/2{\rm\ for\ some\ }i\in W_2\,\}.    \\
\end{align*}
Set
$$
  y_2=\sum_{i\in W_2}\alpha_i x_i.
$$
Then $A_1,A_2,B_1,B_2$, $B_0$, and $C_0$ are mutually disjoint and, as above,
we have
$$
  \V P_{B_1}(y_1+y_2)\V+\V P_{B_2}y_2\V\ge\V P_{B_1}y_1\V-\V P_{A_2}y_2\V. \leqno{(*)}
$$

In general, at the $l^{\rm th}$ step of the induction, we set
\begin{align*}
       W_{l}&=\{\,i\in S_2:x_i(j)=1{\rm\ for\ some\ }j\in B_{l-1}\,\},         \\
       A_{l}&=\{\,j\in A_0:x_i(j)=-1/2{\rm\ for\ some\ }i\in W_l\,\},        \\
\intertext{and}
       B_{l}&=\{\,j\notin A_0:x_i(j)=-1/2{\rm\ for\ some\ }i\in W_l\,\} .   \\
\end{align*}
Set
$$
  y_{l}=\sum_{i\in W_l}\alpha_ix_i.
$$
Then the sets $A_j,B_j$ ($1\le j\le l$), $B_0$, and $C_0$ are mutually disjoint, and we have
$$
  \V P_{B_{l-1}}(y_{l-1}+y_{l})\V+\V P_{B_{l}}y_{l}\V
           \ge\V P_{B_{l-1}}y_{l-1}\V-\V P_{A_{l}}y_{l}\V. \leqno{(*)}
$$

This process must end at some stage $k$ with $W_{k+1}=\emptyset$. Summing the inequalities $(*)$ so obtained
and
simplifying, we have:
\begin{equation}\label{eq: ww}
  \sum_{i=1}^{k} \V P_{B_{i-1}}(y_{i-1}+y_{i})\V+\V P_{B_{k}}y_{k}\V
       \ge\V P_{B_0}y_0\V-\sum_{i=1}^k\V P_{A_{i}}y_{i}\V.
\end{equation}
But, for each $i=0,1,\ldots,k$, we have $B_i\cap{\rm supp}\,y_j =\emptyset$ for
       $j\notin \{\,i,i+1\,\}$. Moreover, the
mutually disjoint sets $A_i$ ($i=1,2,\ldots,k$) are contained in $A_0$. Recalling that $y_0=x$,
from these two observations the inequality \eqref{eq: ww} reduces to
\begin{equation}\label{eq: majorineq}
      \sum_{i=0}^{k}\V P_{B_i}(x+y)\V\ge\V P_{B_0}x\V-\V P_{A_0}y\V.
\end{equation}
Since the
sets $C_0, A_0,B_0,B_1,B_2,\ldots,B_{k}$ are mutually disjoint, we have
$$
  \V x+y\V \ge \V P_{A_0}(x+y) \V +\sum_{i=0}^{k} \V P_{B_i} (x+y) \V +\V P_{C_0}(x+y) \V;
$$
and so, using \eqref{eq: majorineq}:
\begin{equation}\label{eq: MNw}
      \V x+y \V\ge \V P_{A_0}(x+y) \V -\V P_{A_0} y\V +  \V P_{B_0} x\V +\V P_{C_0} x\V.
\end{equation}
Using \eqref{eq: qg} we have $\V P_{A_0}( x+y) \V\ge\V P_{A_0} y\V$, and thus
from \eqref{eq: MNw} we have

\begin{equation}\label{eq: MN}
      \V x+y \V \ge \V P_{B_0} x\V +\V P_{C_0} x\V.
\end{equation}
Also, by \eqref{eq: qg}, we have $2\V P_{A_0} (x+y)\V\ge\V P_{A_0} x\V$; thus,
\begin{equation}\label{eq: L}
      2\V x+y\V\ge\V P_{A_0} x\V.
\end{equation}
Since $\V x\V=\V P_{A_0} x\V+\V P_{B_0} x\V+\V P_{C_0} x\V$,
we may now obtain \eqref{eq: quasigreedy} by adding the inequalities
\eqref{eq: MN} and \eqref{eq: L}.
\end{proof}

\begin{remark}
We note that the sequence of coefficient functionals for $\{\tilde x_i\}_{i=1}^\infty$
is not quasi-greedy. To see this, it is enough to show that
$\{x_i^*\}_{i=1}^\infty$ is not quasi-greedy.
Towards this end, note that
\begin{equation}\label{eq: one}
  \biggl\V\sum_{i=1}^{2^{n+1}-2} (-1)^i x_i^* \biggr\V
    \le \biggl\V\sum_{i=1}^{2^{n+1}-2} (-1)^i y_i^* \biggr\V=1.
\end{equation}
However, defining
$$
  \{\alpha_n\}_{n=1}^\infty=\{\,1,1,\textstyle{1\over2},{1\over2},
  {1\over2},{1\over2},{1\over4},{1\over4},{1\over4},{1\over4},{1\over4},
  {1\over4},{1\over4},{1\over4},\ldots\,\}
$$
and
$$
  z_n=\sum_{i=1}^{2^{n+1}-2} \alpha_i x_i =e_1+e_2+\sum_{i=1}^{2^{n+1}-2}
    {1\over2^{n+1}} e_{2\cdot(2^{n+1}-2)-1+i},
$$
we obtain
\begin{equation}\label{eq: TWo}
      \biggl\V\sum_{i=1}^{2^{n+1}-2} x_i^*\biggr\V\ge{1\over\V z_n\V}
      {|\sum\limits_{i=1}^{2^{n+1}-2}x_i^* (z_n)|}=
      {1\over4}\Bigl<\,\sum\limits_{i=1}^{2^{n+1}-2} x_i^*\,,\, e_1+e_2\Bigr>=\frac{n}{2}.
\end{equation}
It follows from \eqref{eq: one}, \eqref{eq: TWo}, and the Remark of the Introduction
that $\{x_i^*\}_{i=1}^\infty$ is not quasi-greedy.
\end{remark}

\begin{Acknowledgement} The authors thank Professor Wojtaszczyk for providing us
with a preprint of his
paper \cite{W}.
\end{Acknowledgement}


\begin{thebibliography}{99}

\bibitem{HR} J. R. Holub and J. R. Retherford, {\em Some curious bases for $c_0$
and $C[0,1]$\/}, Studia Math. {\bf 34} (1970), 227-240.

\bibitem{L} J. Lindenstrauss, {\em On a certain subspace of $l_1$\/},
     Bull. Acad. Polon. Sci. Ser. sci. math. astr. et phys. {\bf 12} (1964),
     539-542.

\bibitem{LP} J. Lindenstrauss and A. Pe\l czy\'nski, {\em Absolutely summing operators
     in ${\mathcal L}_p$-spaces and their applications\/}, Studia Math.
     {\bf 29} (1968), 275-326.

\bibitem{S} Ivan Singer, {\em Bases in Banach Spaces\/}, vol. I, Springer--Verlag,
    New York-Heidelberg-Berlin, 1970.


\bibitem{W} P. Wojtaszczyk, {\em Greedy algorithm for general biorthogonal systems}, {\em J. Approx. Theory}, to appear.


\end{thebibliography}
\enddocument